%% file: document.tex
\documentclass[proceedings,submission]{dmtcs} 
\usepackage[ascii]{inputenc}
\usepackage[english]{babel}
\usepackage[OT1]{fontenc}
\usepackage{tikz, enumitem, xcolor, shuffle, cancel, extarrows}
\usetikzlibrary{snakes,arrows,positioning,automata,shadows,shapes}
\usepackage[boxed,linesnumbered,ruled,vlined]{algorithm2e}
\DontPrintSemicolon
\SetAlgoCaptionSeparator{:}
\usepackage{rac}

\author[J.-B. \nom{Priez}]{Jean-Baptiste \nom{Priez}}

\title[Lattice of combinatorial Hopf algebras, trees with multiplicities]
{A lattice of combinatorial Hopf algebras,\\
{\LARGE Application to binary trees with multiplicities.}}

\address{\href{mailto:jean-baptiste.priez@lri.fr}{\texttt{jean-baptiste.priez@lri.fr}}\\LRI, B\^at. 650, Universit\'e Paris-Sud 11, 91405 Orsay, France}
\usepackage{setspace}

\begin{document}

\maketitle 
\input{abstract}
\section{Introduction}      \input{intro} 
\section*{Acknowledgements} \input{acknowledgements}
\section{Background}        \input{background/a}

\section{Polynomial realizations and Hopf algebras}\label{sec:rgh} \input{real_good_hopf/intro}
    \subsection{$\phi$-polynomial realization}\input{real_good_hopf/realization}
    \subsection{Alphabet doubling trick}  \input{real_good_hopf/doubl}
    \subsection{Good Hopf algebras}\input{real_good_hopf/good_hopf}

\section{Good monoids}            \input{bon_mono/intro}
    \subsection{Definition}       \input{bon_mono/bons_mono}
    \subsection{Hopf algebra quotient}\input{bon_mono/quotient}
    \subsection{Operations}       \input{bon_mono/prop}

\section{The union of the sylvester and the stalactic congruences}
                      \label{sec:Union} 
                      \input{arbres/intro}
    \subsection{Algorithm and ta\"iga monoid}          \input{arbres/taiga}
    \subsection{Quotient of \wqsym: \pbtm}           \input{arbres/quotient}
    
\section{The hook length formula}      \input{wip/hook}
    
\section{Conclusion, work in progress and perspectives}\input{perspectives}
    
\bibliographystyle{alpha} 
\bibliography{biblio}
\end{document}

%% file: abstract.tex
{ 
\begin{abstract}
    \paragraph{Abstract.}
    In a first part, we formalize the construction of 
    \emph{combinatorial Hopf algebras} from plactic-like monoids using
    \emph{polynomial realizations}. 
    Thank to this construction we reveal a lattice structure on those 
    combinatorial Hopf algebras. As an application, we construct a new 
    combinatorial Hopf algebra on binary trees with multiplicities and use it to 
    prove a \emph{hook length formula} for those trees.

    \textbf{R\'esum\'e.} 
    Dans une premi\`ere partie, nous formalisons la construction 
    d'\emph{alg\`ebres de Hopf combinatoires} \`a partir d'une 
    \emph{r\'ealisation polynomiale} et de mono\"ides de type 
    \emph{mono\"ide plaxique}. Gr\^ace \`a cette construction, nous mettons 
    \`a jour
    une structure de treillis sur ces alg\`ebres de Hopf combinatoires.
    Comme application, nous construisons une nouvelle alg\`ebre de
    Hopf sur des arbres binaires \`a multiplicit\'es et on l'utilise pour 
    d\'emontrer une \emph{formule des \'equerres} sur ces arbres.
\end{abstract}}

%% file: intro.tex
In the past decade a large amount of work in algebraic combinatorics has
been done around combinatorial Hopf algebras. Many have been
constructed on various combinatorial objects such as
\emph{partitions} (symmetric functions \cite{macdonald1995symmetric}), 
\emph{compositions} (\textbf{NCSF} \cite{gelfand1994noncommutative, 
malvenuto1995duality}), 
\emph{permutations} (\fqsym \cite{malvenuto1995duality, 
duchamp2001noncommutative}), 
\emph{set-partitions}  (\wqsym \cite{hivert1999combinatoire,
novelli2005construction}), \emph{binary trees} (\pbt or the \nom{Loday-Ronco} 
Hopf algebra of planar binary trees
\cite{loday1998hopf, hivert2002analogue, hivert2005algebra}), or
\emph{parking functions} (\pqsym 
\cite{novelli2003hopf, novelli2007parking}). 
A powerful method to construct 
those algebras, called \emph{polynomial realization}, is to construct 
the Hopf algebra as a sub algebra of a free algebra of
polynomials (commutative or not) admitting certain
symmetries. Beside the contruction of \emph{Combinatorial Hopf algebra}, 
several recent papers investigate toward the formalization of 
combinatorial applications such as \emph{hook formulas}, or seek 
some structure in this large zoo.

This extended abstract, reports on a work in progress which 
proposes to formalize the construction
of Hopf algebras by polynomial realizations: starting with one of the three 
Hopf algebras \fqsym, \wqsym or \pqsym realized in a \emph{free algebra}, we
impose some relations on the variables. Under some simple
hypotheses, the result is again a Hopf algebra (Theorem 
\ref{thmGoodHopfAlgebra}). Two important examples are already known, 
namely the 
\nom{Poirier-Reutenauer} algebra of tableaux (\fsym 
\cite{poirier1995algebres, duchamp2001noncommutative}) obtained from the 
\emph{plactic monoid} \cite{plaxique} and the \emph{planar binary 
tree algebra} of \nom{Loday-Ronco} obtained from the \emph{sylvester monoid} 
\cite{loday1998hopf, hivert2002analogue, hivert2005algebra}.

We further observe that the construction 
transports the lattice structure on monoids to a lattice structure on those 
Hopf algebras (Theorem \ref{thmBonMonoOpe}). 
This structure was used implicitely by \nom{Giraudo} for 
constructing the Baxter Hopf algebra from the Baxter monoid as the infimum of the 
\emph{sylvester monoid} and its image under \nom{Sch\"utzenberger} 
involution. The supremum of those two monoids is known as the \emph{hypoplactic 
monoid} which gives the algebra of \emph{quasi symmetric functions} 
\cite{novelli2000hypoplactic}.

As an application (Section \ref{sec:Union}) we take the supremum of the
\emph{sylvester monoid} and the \emph{stalactic monoid} of 
\cite{hivert2008commutative}. The result is a monoid on binary 
search trees with multiplicities leading to a Hopf algebra on binary 
trees with multiplicities. Interestingly, there is a \emph{hook length 
formula} for 
those trees (Theorem \ref{thmHookForm}) and we prove it using 
the Hopf algebra as generating series.

%% file: acknowledgements.tex
I would like to thank Florent \nom{Hivert} and Nicolas M.
\nom{Thi\'ery} for their patience and advice during the preparation 
of this manuscript. Several 
examples in this paper were computed using the open-source mathematical 
software \texttt{Sage} \cite{sage} and its extension: 
\texttt{Sage-Combinat} \cite{sagecombinat}. The 
implementation of several Hopf algebras is 
available at \href{http://code.google.com/p/sage-hopf-algebra/}
{\texttt{http://code.google.com/p/sage-hopf-algebra/}} and will be 
available on \texttt{Sage-Combinat} soon and later in \texttt{Sage}.
That code is due to R\'emi \nom{Maurice} and I.

%% file: background/a.tex
In this section, we introduce some notations and three specific
maps from words to words: \emph{standardization}, \emph{packing}, and
\emph{parkization}.
These will be the main tool for polynomial realizations of Hopf algebras.

\subsection{Lattice structure on Congruences} \input{background/cong}

\subsection{Some $\phi$-maps}                 \input{background/appli}

%% file: background/cong.tex
The free monoid $\Ag^*$ on an alphabet $\Ag$ is the set of words
with concatenation as multiplication. We denote by $1$ the empty word.
Recall that a \emph{monoid congruence} is an equivalence relation $\equiv$
 which is left and right compatible with the product;
in other words, for any monoid elements $a,b,c,d$, if $a\equiv b$ and 
$c \equiv d$ then
$a c \equiv b d$.
Starting with two congruences on can build two new congruences:\\[-20pt]
\begin{itemize}[itemsep=-6pt]
  \item the \emph{union}  $\sim \vee \approx$ of $\sim$ and $\approx$ is 
    the transitive closure of the union $\sim$ and $\approx$; that is 
    $u \equiv v$ if there exists $u=u_0, \ldots, u_k =v$
    such that  for any $i$, $u_i \sim u_{i+1}$ or $u_i \approx u_{i+1}$.
    It is the smallest congruence containing both $\sim$ and $\approx$; 
  \item the \emph{intersection} $\sim \wedge \approx$ of $\sim$ and
  $\approx$ is defined as the relation $\equiv$ with $u\equiv v$ if
    $u \sim v$ and $u \approx v$. 
\end{itemize}


%% file: background/appli.tex
Throughout this paper we construct Hopf algebras from the equivalence 
classes of words given by the fibers of some map $\phi$ from
the free monoid to itself. Our main examples are 
\emph{standardization} and \emph{packing} functions which can be defined for any 
totally ordered alphabet $\Ag$. We could easily extend these following 
properties to \emph{parkization} \cite{novelli2003hopf, novelli2007parking}
\emph{if the alphabet $\Ag$ is well-ordered} (any element has a successor).

In the following, \emph{we suppose that $\Ag$ is an totally 
ordered infinite alphabet}. 
Most of the time we use $\Ag = \NN^{>0}$ for simplicity.
For $w$ in $\Ag^*$, we denote by $part(w)$ the ordered set partition 
of positions of $w$ letters obtained as
follows: for each letter $l\in\Ag$ appearing in $w$, there is a 
part containing the positions of each occurrence of $l$ in $w$; the
parts are ordered using the order on the alphabet $\Ag$. 
For example: $part(13231) = [\{1,5\},\{3\},\{2,4\}]$ and $part(1112) =
[\{1,2,3\},\{4\}]$.  

\paragraph{Standardization} $std$ computes the lexicographically smallest 
word $w$ which has same length and same set of inversions. This map is 
used in the realization of the Hopf algebra \fqsym of permutations
\cite{duchamp2001noncommutative, malvenuto1995duality}. 
The image $std(\Ag^*)$ is identified with the set $\mathfrak S$ of all 
permutations.
    
\begin{minipage}[b]{0.607\linewidth} 
    \begin{algorithm}[H]
        \caption{\label{algoStd}Standardization $std$}
        \KwData{$w = (a_1, \ldots, a_k)$ a word of $\Ag^k$}
        \KwResult{$\sigma = (\sigma_1, \ldots, \sigma_k)\in \mathfrak S_k \subseteq \Ag^k$}
 
        $osp \leftarrow part(w)$; $i\leftarrow 1$\;
        \ForAll{$set \in osp$}{
            \ForAll{$p\in set$}{
                $\sigma_p \leftarrow i$\;
                increment($i$)
            }
        }    
       \Return{$\sigma$}
    \end{algorithm}
\end{minipage}
\hfill
\begin{minipage}[b]{0.35\linewidth}
    Some examples:
    
    \begin{tabular}{c|c}
        $w$ & $std(w)$ \\ \hline
        $(7,2,14,3,7)$ & $(3, 1, 5, 2, 4)$\\
        $(23,14,5,92)$ & $(3,2,1,4)$\\
        $(4,2,1,3,5)$ & $(4,2,1,3,5)$\\
        $(1,5,1,1,5,5)$ & $(1, 4, 2, 3, 5, 6)$ 
    \end{tabular}
\end{minipage}
      
\paragraph{Packing} $pack$ computes the lexicographically smallest 
word $w$ which has same ordered set partitions. This map is used in
the realization of the Hopf algebra \wqsym of ordered set partition
or packed words \cite{hivert1999combinatoire,
novelli2005construction}. We identify $tass(\Ag^*)$ with the collection of ordered 
set partitions.
 
\begin{minipage}[b]{0.607\linewidth}
    \begin{algorithm}[H]
       \caption{\label{algoTass}Packing $pack$}
       \KwData{$w = (a_1, \ldots, a_k)$ a word of $\Ag^*$}
       \KwResult{$c = (c_1, \ldots, c_k)$}

       $osp \leftarrow part(w)$; $i \leftarrow 1$\; 
       \ForAll{$set \in osp$}{
            \lForAll{$p\in set$}{
                $c_p \leftarrow i$
            }\;
            increment($i$)
        }    
       \Return{$c$}
    \end{algorithm}
\end{minipage}
\hfill
\begin{minipage}[b]{0.35\linewidth}
    Some examples:
    
    \begin{tabular}{c|c}
        $w$ & $tass(w)$ \\ \hline
        $(3,13,3,2,13)$ & $(2, 3, 2, 1, 3)$\\
        $(2,2,2,5,8,2)$ & $(1, 1, 1, 2, 3, 1)$\\
        $(4,2,1,3,5)$ & $(4,2,1,3,5)$\\
        $(2,3,1,1,2)$ & $(2,3,1,1,2)$
    \end{tabular}
\end{minipage}

Those maps are used to realize some Hopf algebras like 
\fqsym, \wqsym, or \pqsym. For each such map $\phi$ we say that a 
word $w$ is \emph{canonical} if $\phi(w) = w$. For example, 
$1423$ is canonical for $std$ and $1121$ is canonical for $pack$.
The set of canonical words for the \emph{standardization} 
function is the set of permutations set; for the \emph{packing} 
function it is the set of \emph{packed words}. The set $\phi(\Ag^*)$ of 
canonical words is  denoted by $can_\phi$. We call these maps the 
\emph{$\phi$-maps}.

%% file: real_good_hopf/intro.tex
In this section we describe how, from a $\phi$-map, one can construct
a Hopf algebra such as \fqsym, \wqsym, or \pqsym, using two tricks:
\emph{polynomial realization} and \emph{alphabet doubling}.
\emph{Polynomial realizations} are a powerful trick to 
construct algebras as sub-algebras of a free algebra by
manipulating some polynoms having certain symmetries.
Futhermore the \emph{alphabet doubling trick} defines a graded algebra 
morphism on a free algebra which endows it with a compatible coproduct,
that is a Hopf algebra structure. 

%% file: real_good_hopf/realization.tex
The notion of polynomial realizations has been introduced and 
implicitly used in many articles of the ``\textit{phalanst\`ere de 
Marne-la-Vall\'ee}" (France). See \textit{e.g.} \cite{duchamp2001noncommutative, 
novelli605061polynomial, hivert2008commutative}. In the following, we call
\emph{alphabet} $\Ag$ an infinite and totally ordered (when appropriate, we 
assume furthermore that the total order admits a successor function) set of 
symbols all of which are of weight $1$.
By an abuse of language, we call the free algebra the graded algebra 
infinite \textbf{but finite degree} sum of words.


\begin{defi}[Polynomial realization]
    Let $\A := \oplus_{n \geqslant 0} \A_n$ be a graded algebra.
    A \emph{polynomial realization} $r$ of $\A$ is a map which associates
    to each alphabet $\Ag$ an injective graded algebra morphism $r_\Ag$ 
    from $\A$ to the free non-commutative algebra $\KA$ such that, if 
    $\Ag \subset \Bg$, then for all $x\in \A$ one has $r_\Ag (x) = 
    r_\Bg(x)/_\Ag$, where $r_\Bg(x)/_\Ag$ is the sub linear combination 
    obtained from $r_\Bg(x)$ by keeping only those words in $\Ag^*$.
    
    When the realization is clear from
    the context we write $\A(\Ag) := r_\Ag (\A)$ for short. 
\end{defi}

For a given $\phi$, we consider the subspace $\A_\phi$ admitting the basis 
$(m_u)_{u\in can_\phi}$ defined on $\A_\phi (\Ag)$:
\begin{equation}
    \rap[\Ag] (m_u) = \sum_{w\in \Ag^*;\phi(w) = u} w\,.
\end{equation}
The result does not depend on the alphabet.
    For $\phi = std$, $pack$ or $park$ the linear span of 
    $(m_u)_{u\in can_\phi}$
    is a sub-algebra of $\KA$.
\begin{exemple}[Realization of \fqsym]
    If $\phi = std$ then $can_\phi$ is in fact the set of permutations and
     $\Ag_\phi$ is the permutations Hopf algebra \fqsym
    \cite{duchamp2001noncommutative, malvenuto1995duality}. It 
    is realized by the $std$-polynomial realization in $\KA$: 
    let $\GG_\sigma(\Ag):= r_{\Ag,std} (\GG_\sigma)$ such that, for 
    example
    \begin{eqnarray*}\GG_{132}(\NN^*) =  121 + 131 + 132 + 141 + 142 + 143 + \cdots + 
    242 + 243 + \cdots\,.
    \end{eqnarray*}
    
    The realization is an algebra morphism:
    $\GG_\sigma (\Ag) \cdot \GG_\mu(\Ag) = r_{\Ag,std}(\GG_\sigma \times 
    \GG_\mu)$
    where "$\cdot$" is the classical concatenation product on words in
    the free algebra. For example,
    \begin{eqnarray*}
        \mathbb{G}_{213} \times \mathbb{G}_{1} &=& 
\mathbb{G}_{2134} + \mathbb{G}_{2143} + \mathbb{G}_{3142} + \mathbb{G}_{3241}
    \end{eqnarray*}
    which is equivalent to
    \let\boldsymbolOld\boldsymbol
    \renewcommand{\boldsymbol}[1]{\textcolor{red!80}{\boldsymbolOld{#1}}} 
    \begin{eqnarray*}
        r_{std,\NN^*}(\mathbb{G}_{213} \times \mathbb{G}_{1}) &=& 
        \mathbb{G}_{213}(\NN^*) \cdot 
        \boldsymbol{\mathbb{G}_{1}(\NN^*)}\\
        &=& (212 + 213 + 214 + 215 + \cdots) \cdot 
        (\boldsymbol{1} + \boldsymbol{2} + \boldsymbol{3} + \boldsymbol{4} + \boldsymbol{5} + \cdots) \\
        &=& 212\boldsymbol{1} + 212\boldsymbol{2} + 212\boldsymbol{3} + \cdots
        213\boldsymbol{1} + 213\boldsymbol{2} + 213\boldsymbol{3} + \cdots
        + 324\boldsymbol{1} + \cdots
    \end{eqnarray*}
\end{exemple}

\begin{prop}
    If $span((m_u)_{u\in can_\phi})$ is stable under the product 
    $\times$ then it is given by:
    \begin{equation}m_u \times m_v = 
    \sum_{\substack{w:= u'v'\in can_\phi\\ 
            \phi(u') = u; \phi(v')=v}} m_{w}\,.
    \end{equation}
\end{prop}

\begin{remarq}
    \label{rqAlphInfini}
    Let $\Ag, \Bg$ be two totally ordered alphabets such that any
    element in $\Ag$ is strictly smaller than any element of $\Bg$.
    By definition we have the following isomorphisms, where $\sqcup$ 
    denotes the disjoint union:
    \begin{equation}
        \A \simeq \A(\Ag) \simeq \A(\Bg) \simeq \A(\Ag \sqcup \Bg)\,.
    \end{equation}
\end{remarq} 

%% file: real_good_hopf/doubl.tex
The \emph{alphabet doubling trick} \cite{duchamp2001noncommutative, hivert2007introduction}
 is a way to define coproducts.
We consider the algebra $\KA[\Ag\sqcup \Bg]$ generated by two 
{\small(infinite and totally ordered)} alphabets $\Ag$ and $\Bg$ such 
that the letters of $\Ag$ are strictly smaller than the letters of $\Bg$. 
The relation $\commute$ make the letters of $\Ag$ commute with those of $\Bg$. 
One identifies $\KA[\Ag \sqcup \Bg]/_{\commute}$ with the algebra 
$\KA \otimes \KA[\Bg]$.
    We follow here the abuse of language allowing infinite but finite degree 
    sum.
    We denote by $r_{\Ag\sqcup \Bg}(x)/_{\commute}$ the 
    image of $r_{\Ag\sqcup \Bg}(x)$ given by the canonical map from 
    $\KA[\Ag\sqcup \Bg]$ to $\KA[\Ag\sqcup \Bg]/_{\commute}$.
    The map $x \mapsto r_{\Ag \sqcup \Bg} (x)/_{\commute}$ is always
    an algebra morphism from $\A$ to $\KA \otimes \KA[\Bg]$.
    Whenever its image is included in $\A(\Ag)\otimes \A(\Bg)$ this
    defines a coproduct on $\A$.

\begin{defi}[Hopf polynomial realization]
    A \emph{Hopf polynomial realization} $r$ of $\H$ is a polynomial 
    realization such that for all $x$:
    \begin{equation}
        r_{\Ag \sqcup \Bg} (x)/_{\commute} = (r_\Ag \otimes r_\Bg) (\Delta (x))\;. 
    \end{equation}
\end{defi}

\let\boldsymbolOld\boldsymbol
\renewcommand{\boldsymbol}[1]{\textcolor{red!80}{\boldsymbolOld{#1}}}
\begin{exemple}[Coproduct in \fqsym]
    We denote by $\mathbb G_\sigma(\Ag\sqcup\Bg)$ the 
    $std$-polynomial realization of the \fqsym element indexed by 
    $\sigma$ in the algebra $\KA[\Ag\sqcup\Bg]/_{\commute}$. 
    Also we denote by $1,2,3, \cdots$ the symbols of $\Ag$ and in bold red  
    $\boldsymbol{1}, \boldsymbol{2}, \boldsymbol{3}, \boldsymbol{\cdots}$ 
    the symbols of $\Bg$ ordered with $1<2<3<\cdots <\boldsymbol{1} 
    < \boldsymbol{2} < \boldsymbol{3} < \cdots$. Then,
    \begin{eqnarray*}
        \mathbb G_{132} (\Ag \sqcup \Bg) &=&
        121 + 131 + 132 + \cdots + 1\boldsymbol{1} 1 + 1\boldsymbol{1} \boldsymbol{2} + \cdots + 
        1 \boldsymbol{2}\boldsymbol{1} + 1 \boldsymbol{3} \boldsymbol{1} + \cdots 
        + \boldsymbol{1}\boldsymbol{2}\boldsymbol{1} + \cdots\\
        &=& 121 + 131 + 132 + \cdots + 
        11 \cdot \boldsymbol{1} + 11 \cdot \boldsymbol{2} 
        + 12 \cdot \boldsymbol{1} + \cdots
        +\\ &&
        1 \cdot \boldsymbol{2} \boldsymbol{1} + 1\cdot \boldsymbol{3} \boldsymbol{1} + \cdots +
        \boldsymbol{1} \boldsymbol{2} \boldsymbol{3} + \boldsymbol{1} \boldsymbol{3}\boldsymbol{2} +
        \boldsymbol{2}\boldsymbol{3}\boldsymbol{2} + \cdots \\
        &=& \Delta(\GG_{132}) = 1 \otimes \mathbb{G}_{132} + \mathbb{G}_{1} 
        \otimes \mathbb{G}_{21} + \mathbb{G}_{12} \otimes \mathbb{G}_{1} + 
        \mathbb{G}_{132} \otimes 1\;.
    \end{eqnarray*}
\end{exemple}

%% file: real_good_hopf/good_hopf.tex
We call a \emph{Hopf algebra $\H_\phi$} \emph{good} if it is defined
by a Hopf polynomial realization $r_\phi$. 
We call a \emph{function $\phi$} \emph{good}
if it produces a \emph{good Hopf algebra} $\H_\phi$. 
Currently, we know three main good Hopf algebras: \fqsym, \wqsym 
and \pqsym are respectivly associated to the \emph{standardization},
\emph{packing} and \emph{parkization} functions.

%% file: bon_mono/intro.tex
In the previous section (Section \ref{sec:rgh}), we realized some 
Hopf algebras in free algebras. In this section, we give sufficient
conditions on a congruence $\equiv$ to build a combinatorial quotient of a 
\emph{good Hopf algebra}. We call a \emph{monoid good} if it statisfies these 
conditions. 
The first condition is about the $\phi$-map used to realize the 
\emph{good Hopf algebra} in free algebras. We give a sufficient 
compatibility between $\phi$ and $\equiv$ to ensure the product is carried to 
the quotient.
The second condition ensures that the \emph{alphabet doubling trick} map.
It is used to project the coproduct in the quotient. 
Under these conditions, a monoids is guaranted to produce a Hopf algebra 
quotient (Theorem \ref{thmGoodHopfAlgebra}).   
Furthermore, these conditions on monoid are preserved under taking
infimum and supremum (Theorem \ref{thmBonMonoOpe}).

%% file: bon_mono/bons_mono.tex
The notion of \emph{Good monoids} has been introduced by \nom{Hivert-Nzeutchap}
\cite{hivert5combinatoire, hivertnzeutchap2007} to build quotients 
(sub-algebras) of \fqsym. We could also mention PhD thesis \cite{giraudo}. 
 
A \emph{good monoid} is a monoid which has similar properties, as 
the \emph{plactic monoid} \cite{plaxique, knuth1970permutations}.
We consider a free monoid $\Ag^*$ with concatenation product "$\cdot$",
a congruence $\equiv$ on $\Ag^*$ and a map $\phi: \Ag^* \to \Ag^*$.
We define the \emph{evaluation} $ev(w)$ of a word $w$ as its number of
occurrences of each letter of $w$. For example, the words \texttt{ejajv} 
and \texttt{jjaev} have the same evaluation: both have one \texttt{a}, one \texttt{e}, 
one \texttt{v} and two \texttt{j}.
The free monoid $\Ag^*/\equiv$ is a \emph{$\phi$-good monoid} if it has 
the following properties:
\begin{defi}[$\phi$-congruence] 
    The congruence $\equiv$ is a \emph{$\phi$-congruence} if for all
    $u,v\in\Ag^*$,  
	$u \equiv v$ if and only if $\phi (u) \equiv \phi(v)$ and 
	$ev(u) = ev(v)$.
\end{defi}
This first compatibility is sufficient to build a quotient algebra of
$\A_\phi$.

\begin{defi}[Compatibility with restriction to alphabet intervals]
    The congruence $\equiv$ is compatible with the \emph{restriction to 
    alphabet intervals} if, for all $u,v\in\Ag^*$ such that $u\equiv v$ one has  
    $u_{|I} \equiv v_{|I}$ for any $I$ interval of $\Ag$, where 
    $w_{|\Ag}$ is word restricted to the alphabet $\Ag$.
\end{defi}

This second compatibility in association with the first ensures that
\emph{alphabet doubling trick} defines a quotient coproduct.
Both compatibilities give us an extended definition of a 
\nom{Hivert-Nzeutchap}'s \emph{good monoid} which one is defined only with $\phi$
the standardization map:
\begin{defi}[$\phi$-good monoid]
    A quotient $\Ag^*/\equiv$ of the free monoid is a 
    \emph{$\phi$-good monoid} 
    if $\equiv$ is a $\phi$-congruence and is compatible with 
    restriction to alphabet intervals. 
    We call such a congruence a \emph{$\phi$-good congruence}.
\end{defi}

In the following examples, we denote words 
of $\Ag^*$ by $u,v,w$ and the letters by $a,b,c$.  

\begin{exemple}[sylvester and stalactic monoids]
    \label{exBonMono}
    The sylvester congruence: $\equiv_{sylv}$, defined by
    \begin{eqnarray}
        u \cdot ac \cdot w \cdot b \cdot v \equiv_{sylv}
        u \cdot ca \cdot w \cdot b \cdot v \text{ whenever }
        a \leqslant b < c\;, \label{eqnSyl}
    \end{eqnarray}
    is $std$-compatible and compatible with the restriction to alphabet intervals.
    Thanks to the binary search tree insertion algorithm the equivalence classes 
    are in natural bijection with binary search trees.
    The quotient monoid is a monoid on binary search trees called the
    sylvester monoid in \cite{hivert2002analogue, hivert2005algebra}.
    
    The stalactic congruence \cite{hivert2008commutative}: $\equiv_{stal}$,
    defined by
    \begin{eqnarray}
        u\cdot ba \cdot v \cdot b \cdot w \equiv_{stal}
        u\cdot ab \cdot v \cdot b \cdot w\;, \label{eqnStal}
    \end{eqnarray}
    is compatible with \emph{packing} but not with
    \emph{standardization}. The quotient monoid is the \emph{stalactic 
    monoid}. 
    It is clear that any stalactic
    class contains a word of the form $a_1^{m_1}a_2^{m_2}\ldots a_k^{m_k}$,
    where the $a_i$ are distinct. We call these words 
    canonical. We represent a stalactic class with a planar diagram 
    such that, in any column, the boxes contain the same letter.
    \begin{center}
        $51543151145312455\equiv_{stal} 3^2 1^5 2^1 4^3 5^6 \longleftrightarrow \vcenter{\hbox{\begin{tikzpicture}[scale=.3]
            \foreach \x in {0,-1} {
                \draw (0,\x) rectangle (1,\x + 1);
                \draw (.5,\x+.5) node{$\scriptstyle 3$};
            }
            \foreach \x in {0,-1,-2,-3,-4} {
                \draw (1,\x) rectangle (2,\x + 1);
                \draw (1.5,\x+.5) node{$\scriptstyle 1$};
            }
            \foreach \x in {0} {
                \draw (2,\x) rectangle (3,\x + 1);
                \draw (2.5,\x+.5) node{$\scriptstyle 2$};
            }
            \foreach \x in {0,-1,-2} {
                \draw (3,\x) rectangle (4,\x + 1);
                \draw (3.5,\x+.5) node{$\scriptstyle 4$};
            }
            \foreach \x in {0,-1,-2,-3,-4,-5} {
                \draw (4,\x) rectangle (5,\x + 1);
                \draw (4.5,\x+.5) node{$\scriptstyle 5$};
            }
        \end{tikzpicture}}}$
    \end{center}  
\end{exemple}

%% file: bon_mono/quotient.tex
These differents \emph{good monoids} tools was used to 
(re-)define several Hopf algebra quotients: \fsym 
the Free Symmetric functions Hopf algebra \cite{duchamp2001noncommutative}, 
\pbt \cite{loday1998hopf, hivert2002analogue, hivert2005algebra} or 
Baxter Hopf algebra 
\cite{giraudo2011algebraic, giraudo2012algebraic}; the Hopf algebra 
associated with the stalactic monoid \cite{hivert2008commutative};
or \cqsym \cite{novelli2003hopf, novelli2007parking} (a \pqsym quotient).\par

\begin{lemme}[Algebra quotient]
    \label{lemmeAlgebraQuotient}
    Let $\H_\phi$ be a \emph{good Hopf algebra} and $\equiv$ be a 
    $\phi$-good congruence such that its free monoid quotient is a 
    \emph{$\phi$-good monoid}.
    Then, the quotient $\H_\phi/_\equiv$ is an algebra quotient whose
    bases are indexed by $can_\phi/_\equiv$, identifying basis elements 
    $m_u$ and $m_v$ whenever $u \equiv v$.
\end{lemme}

\newcommand{\QQm}[1]{\QQ_{\vcenter{\hbox{\input{bon_mono/exStal/#1}}}}}
\begin{exemple}[\pbt and Hopf algebra stalactic]
    We go back to Example \ref{exBonMono}.
    The sylvester quotient of \fqsym is the Hopf algebra \pbt 
    \cite{loday1998hopf, hivert2002analogue, hivert2005algebra}.
    
    The stalactic monoid gives a quotient of \wqsym. Let $\pi$ be the
    projection of \wqsym
    in $\wqsym/_{\equiv_{stal}}$ and $u := 112$ and $v := 11$ two 
    (packed) words. We denote by $\pi$ the projection of $\MM_u$ by $\QQ_s$,
    with $s$ the planar diagram associated to the stalactic class of 
    $u$. 
    \begin{eqnarray*}
        \pi(\mathbb{M}_{112} \times \mathbb{M}_{11}) &=& 
        \pi(\mathbb{M}_{11211} + \mathbb{M}_{11222} + \mathbb{M}_{11233} + 
        \mathbb{M}_{11322} + \mathbb{M}_{22311})\\
        = \ensuremath{\QQ^m_{\vcenter{\hbox{\scalebox{.4}{\input{arbres/quotient/prod/a1}}}}}} \times \ensuremath{\QQ^m_{\vcenter{\hbox{\scalebox{.4}{\input{arbres/quotient/prod/a2}}}}}} &=& 
        \;\;\;\ensuremath{\QQ^m_{\vcenter{\hbox{\scalebox{.4}{\input{arbres/quotient/prod/e1}}}}}} + \ensuremath{\QQ^m_{\vcenter{\hbox{\scalebox{.4}{\input{arbres/quotient/prod/e2}}}}}} + \ensuremath{\QQ^m_{\vcenter{\hbox{\scalebox{.4}{\input{arbres/quotient/prod/e3}}}}}} + \ensuremath{\QQ^m_{\vcenter{\hbox{\scalebox{.4}{\input{arbres/quotient/prod/e4}}}}}} + \ensuremath{\QQ^m_{\vcenter{\hbox{\scalebox{.4}{\input{arbres/quotient/prod/e5}}}}}}  
    \end{eqnarray*}
\end{exemple}

\begin{lemme}[Coalgebra quotient]
    \label{lemmeCogebraQuotient}
    The quotient $\H/_\equiv$ is a coalgebra quotient.
\end{lemme}

\begin{preuve}
    The relation $\equiv$ is compatible with the restriction to 
    alphabet intervals, hence the alphabet doubling trick ensures that
    coproduct projects to the quotient.
\end{preuve}

\begin{exemple}
    \begin{eqnarray*}
        \pi(\Delta(\mathbb{M}_{332122})) &=& \pi(
        1 \otimes \mathbb{M}_{332122} + \mathbb{M}_{1} \otimes 
        \mathbb{M}_{22111} + \mathbb{M}_{2122} \otimes \mathbb{M}_{11} + 
        \mathbb{M}_{332122} \otimes 1)\\
        = \Delta(\ensuremath{\QQ^m_{\vcenter{\hbox{\scalebox{.4}{\input{arbres/quotient/prod/c1}}}}}}) &=& 1 \otimes \ensuremath{\QQ^m_{\vcenter{\hbox{\scalebox{.4}{\input{arbres/quotient/prod/c1}}}}}} +
        \ensuremath{\QQ^m_{\vcenter{\hbox{\scalebox{.4}{\input{arbres/quotient/prod/c2}}}}}} \otimes \ensuremath{\QQ^m_{\vcenter{\hbox{\scalebox{.4}{\input{arbres/quotient/prod/c3}}}}}} + \ensuremath{\QQ^m_{\vcenter{\hbox{\scalebox{.4}{\input{arbres/quotient/prod/c4}}}}}} \otimes \ensuremath{\QQ^m_{\vcenter{\hbox{\scalebox{.4}{\input{arbres/quotient/prod/c5}}}}}} + 
        \ensuremath{\QQ^m_{\vcenter{\hbox{\scalebox{.4}{\input{arbres/quotient/prod/c1}}}}}} \otimes 1
    \end{eqnarray*}
\end{exemple}

\begin{thm}[Good monoid and good Hopf algebra]
    \label{thmGoodHopfAlgebra}
    Let $\H_\phi$ be a good Hopf algebra and $\equiv$ be a $\phi$-good
    congruence. 
    The quotient $\H/_\equiv$ is a Hopf algebra quotient.
\end{thm}

\begin{coro}
    The dual Hopf algebra $(\H/_{\equiv})^\#$ is a sub-algebra of
    the dual Hopf algebra $\H^\#$, with basis given by:
    \begin{equation}
  \bar M^\#_{U\in can_\phi/_{\equiv}} = \sum_{u\in U} m_u^\#\;.
    \end{equation} 
\end{coro}

%% file: arbres/quotient/prod/a1.tex
{\newcommand{\AOnooo}{\node (0) [shape=rectangle, rounded corners, draw=blue!80, ultra thick] 
{$1$};}\newcommand{\AOnooa}{\node (1) [shape=rectangle, draw=orange!80, ultra thick] 
{$2$};}\newcommand{\AOnoob}{\node (2) [shape=ellipse, draw=c1,ultra thick, decorate, decoration=zigzag] 
{$1$};}\begin{tikzpicture}[auto]
\matrix[column sep=.1cm, row sep=.1cm]{
         & \AOnooo \\ 
 \AOnooa &         & \AOnoob \\ 
};
\path[red, ultra thick] (0) edge (2) edge (1);
\end{tikzpicture}}

%% file: arbres/quotient/prod/a2.tex
{\newcommand{\AOnooo}{\node (0) [circle, draw=c3, ultra thick] 
{$1$};}\begin{tikzpicture}[auto]
\matrix[column sep=.1cm, row sep=.1cm]{ \AOnooo \\  };
\end{tikzpicture}}

%% file: bon_mono/exStal/e1.tex
\begin{tikzpicture}[scale=.3]
    \foreach \x in {0} {
        \draw (0,\x) rectangle (1,\x + 1);
        \draw (.5,\x+.5) node{$\scriptstyle 2$};
    }
    \foreach \x in {0,-1,-2,-3} {
        \draw (1,\x) rectangle (2,\x + 1);
        \draw (1.5,\x+.5) node{$\scriptstyle 1$};
    }
\end{tikzpicture}

%% file: bon_mono/exStal/e2.tex
\begin{tikzpicture}[scale=.3]
    \foreach \x in {0,-1} {
        \draw (0,\x) rectangle (1,\x + 1);
        \draw (.5,\x+.5) node{$\scriptstyle 1$};
    }
    \foreach \x in {0,-1,-2} {
        \draw (1,\x) rectangle (2,\x + 1);
        \draw (1.5,\x+.5) node{$\scriptstyle 2$};
    }
\end{tikzpicture}

%% file: bon_mono/exStal/e3.tex
\begin{tikzpicture}[scale=.3]
    \foreach \x in {0,-1} {
        \draw (0,\x) rectangle (1,\x + 1);
        \draw (.5,\x+.5) node{$\scriptstyle 1$};
    }
    \foreach \x in {0} {
        \draw (1,\x) rectangle (2,\x + 1);
        \draw (1.5,\x+.5) node{$\scriptstyle 2$};
    }
    \foreach \x in {0,-1} {
        \draw (2,\x) rectangle (3,\x + 1);
        \draw (2.5,\x+.5) node{$\scriptstyle 3$};
    }
\end{tikzpicture}

%% file: bon_mono/exStal/e4.tex
\begin{tikzpicture}[scale=.3]
    \foreach \x in {0,-1} {
        \draw (0,\x) rectangle (1,\x + 1);
        \draw (.5,\x+.5) node{$\scriptstyle 1$};
    }
    \foreach \x in {0} {
        \draw (1,\x) rectangle (2,\x + 1);
        \draw (1.5,\x+.5) node{$\scriptstyle 3$};
    }
    \foreach \x in {0,-1} {
        \draw (2,\x) rectangle (3,\x + 1);
        \draw (2.5,\x+.5) node{$\scriptstyle 2$};
    }
\end{tikzpicture}

%% file: bon_mono/exStal/e5.tex
\begin{tikzpicture}[scale=.3]
    \foreach \x in {0,-1} {
        \draw (0,\x) rectangle (1,\x + 1);
        \draw (.5,\x+.5) node{$\scriptstyle 2$};
    }
    \foreach \x in {0} {
        \draw (1,\x) rectangle (2,\x + 1);
        \draw (1.5,\x+.5) node{$\scriptstyle 3$};
    }
    \foreach \x in {0,-1} {
        \draw (2,\x) rectangle (3,\x + 1);
        \draw (2.5,\x+.5) node{$\scriptstyle 1$};
    }
\end{tikzpicture}

%% file: bon_mono/exStal/c1.tex
\begin{tikzpicture}[scale=.3]
    \foreach \x in {0,-1} {
        \draw (0,\x) rectangle (1,\x + 1);
        \draw (.5,\x+.5) node{$\scriptstyle 3$};
    }
    \foreach \x in {0} {
        \draw (1,\x) rectangle (2,\x + 1);
        \draw (1.5,\x+.5) node{$\scriptstyle 1$};
    }
    \foreach \x in {0,-1,-2} {
        \draw (2,\x) rectangle (3,\x + 1);
        \draw (2.5,\x+.5) node{$\scriptstyle 2$};
    }
\end{tikzpicture}

%% file: bon_mono/exStal/c2.tex
\begin{tikzpicture}[scale=.3]
    \foreach \x in {0} {
        \draw (0,\x) rectangle (1,\x + 1);
        \draw (.5,\x+.5) node{$\scriptstyle 1$};
    }
\end{tikzpicture}

%% file: bon_mono/exStal/c3.tex
\begin{tikzpicture}[scale=.3]
    \foreach \x in {0,-1} {
        \draw (0,\x) rectangle (1,\x + 1);
        \draw (.5,\x+.5) node{$\scriptstyle 2$};
    }
    \foreach \x in {0,-1,-2} {
        \draw (1,\x) rectangle (2,\x + 1);
        \draw (1.5,\x+.5) node{$\scriptstyle 1$};
    }
\end{tikzpicture}

%% file: bon_mono/exStal/c4.tex
\begin{tikzpicture}[scale=.3]
    \foreach \x in {0} {
        \draw (0,\x) rectangle (1,\x + 1);
        \draw (.5,\x+.5) node{$\scriptstyle 1$};
    }
    \foreach \x in {0,-1,-2} {
        \draw (1,\x) rectangle (2,\x + 1);
        \draw (1.5,\x+.5) node{$\scriptstyle 2$};
    }
\end{tikzpicture}

%% file: bon_mono/exStal/c5.tex
\begin{tikzpicture}[scale=.3]
    \foreach \x in {0,-1} {
        \draw (0,\x) rectangle (1,\x + 1);
        \draw (.5,\x+.5) node{$\scriptstyle 1$};
    }
\end{tikzpicture}

%% file: bon_mono/prop.tex
Previously we introduced some \emph{good functions} $\phi$: $std$, 
$pack$ (and $park$). It is interesting to investigate the connections 
between them:

\begin{defi}[refinement]
    Let $\phi$ and $\pi$ be two functions.
    We say that $\pi$ refines $\phi$, written $\phi\prec \pi$ 
    if $\phi(\pi(u)) = \phi(u)$ for all $u\in \Ag^*$.
\end{defi}

It is clear that \emph{refinement} is an order.

\begin{prop}[$std$, $tass$, $park$ and refinement]
    \label{propStdTassParkInclusion}
    For these three functions: \emph{standardization} $std$, 
    \emph{packing} $pack$ and \emph{parking} $park$ we have the
    relation:
      $std \prec pack \prec park$. 
\end{prop}

\begin{prop}[Good functions and refinement]
    \label{propBonneFoncInclusion}
    Let $\phi$ and $\pi$ be two \emph{good functions} such that
    $\phi \prec \pi$.
    Then any \emph{$\phi$-good monoid} is a \emph{$\pi$-good monoid}.
\end{prop}

Propositions \ref{propStdTassParkInclusion} and
\ref{propBonneFoncInclusion} give us, for example, that any $std$-good 
monoid is $pack$-good.
Furthermore operations on two \emph{good congruences} give \emph{good
congruences}.

\begin{thm}[$\vee$, $\wedge$ and good congruences]
    \label{thmBonMonoOpe}
    The union and intersection of two \emph{$\phi$-good 
    congruences} $\sim$ and $\approx$ are \emph{$\phi$-good congruences}.
\end{thm}

As an intriguing consequence the lattice structure on monoids is 
transported to Hopf algebras. Several examples of this are know.

\begin{exemple}
    The intersection ($\equiv_{sylv} \wedge \equiv_{\#sylv}$) of 
    the sylvester relation (\ref{eqnSyl}) and its image under the 
    \nom{Sch\"utzenberger} involution gives $std$-good monoid: the \emph{Baxter monoid}
    \cite{giraudo2011algebraic, giraudo}.
    
    The union ($\equiv_{sylv} \vee \equiv_{\#sylv}$) of those relations
    gives the \emph{hypoplactic monoid} \cite{novelli2000hypoplactic}.
\end{exemple}

In the sequel, we study in detail another example.

%% file: arbres/intro.tex
As an application of the preceding construction, we consider the union 
($\equiv_{sylv} \vee \equiv_{stal}$) of the sylvester congruence (\ref{eqnSyl}) 
and the stalactic congruence (\ref{eqnStal}); we call it
the \emph{ta\"iga relation} $\equiv_t$,
    \begin{eqnarray}
    \label{eqnTaig}\begin{split}
        u \cdot ac \cdot v \cdot b \cdot w \equiv_{t}&\ 
        u \cdot ca \cdot v \cdot b \cdot w &\text{ for }
        a \leqslant b < c\;, \\
        u\cdot ba \cdot v \cdot b \cdot w \equiv_{t}&\ 
        u\cdot ab \cdot v \cdot b \cdot w&
    \end{split}\end{eqnarray}

From Proposition \ref{propBonneFoncInclusion} we know that the sylvester 
congruence (\ref{eqnSyl}) is a \emph{$pack$-good congruence} and from 
Theorem \ref{thmBonMonoOpe} we deduce
that the \emph{ta\"iga monoid} is a \emph{$pack$-good monoid}.

\begin{figure}
  \centering
      $w:=\texttt{45142234212}\;\;\left\{\begin{matrix}
        \P(w) = \;\; & \vcenter{\hbox{\scalebox{.6}{\input{arbres/fig1}}}} \;\; = \B(w)\\
        \Q(w) = \;\; & \vcenter{\hbox{\scalebox{.6}{\input{arbres/fig2}}}}
      \end{matrix}\right.$
  \caption{\label{figABRMa}We start by considering the packed word 
  \texttt{45142234212}, and insert it in a BSTM by the algorithm $\P$; that give us 
    $\P(\texttt{45142234212})$ above in the middle. On the right, there is 
    the BTm ($\B(w)$) associated with the BSTM ($\P(w)$) of 
    $\wqsym/_{\equiv_t}$. At the top of the figure there is the $P$-symbol
    given by $\P$ or $\B$ and below the $Q$-symbol is given by $\Q$.}
\end{figure}
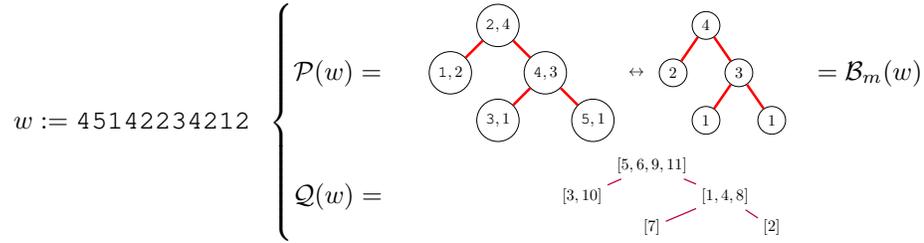 

%% file: arbres/fig1.tex
{\newcommand{\AOnooo}{\node (0) [circle, fill=white, text=black, draw] 
{$\texttt{2},4$};}\newcommand{\AOnooa}{\node (1) [circle, fill=white, text=black, draw] 
{$\texttt{1},2$};}\newcommand{\AOnoob}{\node (2) [circle, fill=white, text=black, draw] 
{$\texttt{3},1$};}\newcommand{\AOnooc}{\node (3) [circle, fill=white, text=black, draw] 
{$\texttt{4},3$};}\newcommand{\AOnood}{\node (4) [circle, fill=white, text=black, draw] 
{$\texttt{5},1$};}\newcommand{\AOnooe}{\node (5) [circle, fill=white, text=black, draw] 
{$2$};}\newcommand{\AOnoof}{\node (6) [circle, fill=white, text=black, draw] 
{$4$};}\newcommand{\AOnoog}{\node (7) [circle, fill=white, text=black, draw] 
{$3$};}\newcommand{\AOnooh}{\node (8) [circle, fill=white, text=black, draw] 
{$1$};}\newcommand{\AOnooi}{\node (9) [circle, fill=white, text=black, draw] 
{$1$};}\begin{tikzpicture}[auto]
\matrix[column sep=.1cm, row sep=.1cm]{
                & \AOnooo &                &               & &               & &                & \AOnoof \\ 
 \AOnooa &                & \AOnooc &                &     &  \node{$\leftrightarrow$};   &      & \AOnooe &                & \AOnoog \\ 
                & \AOnoob &                & \AOnood &    &            &       &         & \AOnooh &                & \AOnooi \\ 
};
\path[red, ultra thick] (7) edge node {} (9);
\path[red, ultra thick] (3) edge node {} (2);
\path[red, ultra thick] (0) edge node {} (1);
\path[red, ultra thick] (3) edge node {} (4);
\path[red, ultra thick] (7) edge node {} (8);
\path[red, ultra thick] (0) edge node {} (3);
\path[red, ultra thick] (6) edge node {} (7);
\path[red, ultra thick] (6) edge node {} (5);
\end{tikzpicture}}

%% file: arbres/fig2.tex
{\newcommand{\AOnooo}{\node (0)  
{$[5,6,9,11]$};}\newcommand{\AOnooa}{\node (1) 
{$[3,10]$};}\newcommand{\AOnoob}{\node (2) 
{$[7]$};}\newcommand{\AOnooc}{\node (3)  
{$[1,4,8]$};}\newcommand{\AOnood}{\node (4) 
{$[2]$};}\begin{tikzpicture}[auto]
\matrix[column sep=.1cm, row sep=.1cm]{
                & \AOnooo &                &              \\ 
 \AOnooa &                & \AOnooc &                \\ 
                & \AOnoob &                & \AOnood \\ 
};
\path[purple,thick] (3) edge node {} (2);
\path[purple,thick] (0) edge node {} (1);
\path[purple,thick] (3) edge node {} (4);
\path[purple,thick] (0) edge node {} (3);
\end{tikzpicture}}

%% file: arbres/taiga.tex
The ta\"iga congruence can be calculated using an insertion algorithm
similar to the binary search tree insertion 
(see Algorithm \ref{algoInsertion} for a definition). 
This insertion algorithm uses a search tree structure:

\begin{defi}[Binary search tree with multiplicity]
    A (planar) \emph{binary search tree with multiplicity} (BS\-TM) is a binary  
    tree $T$ where each node is labelled by a letter $l$ and a non-negative 
    integer $k$, called 
    the multiplicity, so that $T$ is a binary search tree if we drop the
    multiplicities and such that each letter appears  at most once in $T$.
    
    We denote by $(l,k)$ a node label and for any node $n$, by $l(n)$ its 
    letter and by $m(n)$ its multiplicity.   
\end{defi}

\begin{figure}
\begin{minipage}[b]{.50\linewidth} 
\begin{algorithm}[H]
   \caption{\label{algoInsertion}insertion in a BSTM}
   \KwData{$t$ a BSTM with $L_t$ and $R_t$ its left and right subtrees,
        and $l$ a letter of $\Ag$}
   \KwResult{$t$ with $l$ inserted}
 
   \BlankLine
   \eIf{$t$ is empty tree}{
        $t \leftarrow$ node labelled by $(l,1)$
   }{  \eIf{$l(t) = l$}{
            increment $m(t)$\;
       }{
            \lIf{$l(t)<l$}{ insert recursively $l$ in $L_t$ }\;
            \lElse{ insert recursively $l$ in $R_t$ }
   }   }
   \Return{$t$}
\end{algorithm}
\end{minipage}\hfill
\begin{minipage}[b]{.42\linewidth}
Insertion of word \texttt{541214} from \emph{the right to the left} in 
the empty tree :
    
\begin{tabular}{lclc}
    $\xlongrightarrow{\texttt{4}}$ & $\vcenter{\hbox{\scalebox{.6}{\begin{tikzpicture}[auto]
        \node[circle,draw=orange!80,ultra thick]{$\texttt{4},1$};
    \end{tikzpicture}}}}$    &
    $\xlongrightarrow{\texttt{1}}$ & $\vcenter{\hbox{\scalebox{.6}{\input{arbres/d/a1}}}}$\\ \hline
    $\xlongrightarrow{\texttt{2}}$ & $\vcenter{\hbox{\scalebox{.6}{\input{arbres/d/a2}}}}$ &
    $\xlongrightarrow{\texttt{1}}$ & $\vcenter{\hbox{\scalebox{.6}{\input{arbres/d/a3}}}}$\\ \hline
    $\xlongrightarrow{\texttt{4}}$ & $\vcenter{\hbox{\scalebox{.6}{\input{arbres/d/a4}}}}$ &
    $\xlongrightarrow{\texttt{5}}$ & $\vcenter{\hbox{\scalebox{.6}{\input{arbres/d/a5}}}}$
\end{tabular}
\end{minipage}
\end{figure}

We denote by $\P(w)$ the result of the insertion using Algorithm 
\ref{algoInsertion} of $w$ from \emph{the right to the left} in the 
empty tree (\textit{cf}. the left part of the figure \ref{figABRMa}).

\begin{prop}
    \label{propEquivSsiP}
    The ta\"iga classes are the fibers of $\P$. That is for $u$ 
    and $v$ two words: 
        $u \equiv_t v$ if and only if $\P(u) = \P(v)$.
\end{prop}

The $Q$-symbol of $w$ is the tree $\Q(w)$ of same shape as $\P(w)$ which 
records the 
positions of each inserted letter. This gives us a \emph{\nom{Robinson-Schensted} like 
correspondance} \cite{lothaire2002algebraic} (\textit{cf}. Figure 
\ref{figABRMa}). As a corollary of Theorem \ref{thmBonMonoOpe} we get

\begin{coro}
    The ta\"iga monoid is a \emph{$tass$-good monoid}.
\end{coro}

%% file: arbres/d/a1.tex
\newcommand{\AOnooo}{\node (0) [circle, fill=white, text=black, draw] {$\texttt{4},1$};}\newcommand{\AOnooa}{\node (1) [circle,draw=orange!80,ultra thick] {$\texttt{1},1$};}\begin{tikzpicture}[auto]
\matrix[column sep=.1cm, row sep=.1cm]{
                & \AOnooo \\ 
 \AOnooa \\ 
};
\path[ultra thick, red] (0) edge node {} (1);
    \end{tikzpicture}

%% file: arbres/d/a2.tex
\newcommand{\AOnooo}{\node (0) [circle, fill=white, text=black, draw] {$\texttt{4},1$};}\newcommand{\AOnooa}{\node (1) [circle, fill=white, text=black, draw] {$\texttt{1},1$};}\newcommand{\AOnoob}{\node (2) [circle, draw=orange!80,ultra thick] {$\texttt{2},1$};}\begin{tikzpicture}[auto]
\matrix[column sep=.1cm, row sep=.1cm]{
                & \AOnooo \\ 
 \AOnooa \\
                & \AOnoob \\ 
};
\path[ultra thick, red] (0) edge (1)
                        (1) edge (2);
    \end{tikzpicture}

%% file: arbres/d/a3.tex
\newcommand{\AOnooo}{\node (0) [circle, fill=white, text=black, draw] {$\texttt{4},1$};}\newcommand{\AOnooa}{\node (1) [circle, draw=orange!80,ultra thick] {$\texttt{1},2$};}\newcommand{\AOnoob}{\node (2) [circle, fill=white, text=black, draw] {$\texttt{2},1$};}\begin{tikzpicture}[auto]
\matrix[column sep=.1cm, row sep=.1cm]{
                & \AOnooo \\ 
 \AOnooa \\
                & \AOnoob \\ 
};
\path[ultra thick, red] (0) edge (1)
                        (1) edge (2);
    \end{tikzpicture}

%% file: arbres/d/a4.tex
\newcommand{\AOnooo}{\node (0) [circle, draw=orange!80,ultra thick] {$\texttt{4},2$};}\newcommand{\AOnooa}{\node (1) [circle, fill=white, text=black, draw] {$\texttt{1},2$};}\newcommand{\AOnoob}{\node (2) [circle, fill=white, text=black, draw] {$\texttt{2},1$};}\begin{tikzpicture}[auto]
\matrix[column sep=.1cm, row sep=.1cm]{
                & \AOnooo \\ 
 \AOnooa \\
                & \AOnoob \\ 
};
\path[ultra thick, red] (0) edge (1)
                        (1) edge (2);
    \end{tikzpicture}

%% file: arbres/d/a5.tex
\newcommand{\AOnooo}{\node (0) [circle, fill=white, text=black, draw] {$\texttt{4},2$};}\newcommand{\AOnooa}{\node (1) [circle, fill=white, text=black, draw] {$\texttt{1},2$};}\newcommand{\AOnoob}{\node (2) [circle, fill=white, text=black, draw] {$\texttt{2},1$};}\newcommand{\AOnooc}{\node (3) [circle, draw=orange!80,ultra thick] {$\texttt{5},1$};}\begin{tikzpicture}[auto]
\matrix[column sep=.1cm, row sep=.1cm]{
                & \AOnooo \\ 
 \AOnooa        &             & \AOnooc \\
                & \AOnoob \\ 
};
\path[ultra thick, red] (0) edge (1) edge (3)
                        (1) edge (2);
    \end{tikzpicture}

%% file: arbres/quotient.tex
Thanks to Proposition \ref{propEquivSsiP}, the set of packed words giving 
the same tree by algorithm $\P$ is exactly a ta\"iga class of packed 
words. As in \cite{hivert2002analogue}, we consider a binary trees with
multiplicities without letters.

\begin{defi}[BTM]
    A \emph{binary tree with multiplicities} (BTM) is a (planar) binary 
    tree labelled by non-negative integers on its nodes.
    The size of a BTM $T$ denoted by $\vabs T$ is the sum of the 
    multiplicities. 
\end{defi}

Let $T_w$ be a BSTM associated to a packed word $w$, and $T$ be
the BTM obtained by removing its letters. One can recover uniquely 
$T_w$ from $T$: indeed each letter of $T_w$
is deduced by a left infix reading of $T$.
We identify the set of words in $pack(\Ag^*)/_{\equiv_t}$ of size $k$ 
(for $k\geqslant 0$) with the set of BTM of size $k$. 
We denote $\B$ the algorithm which 
computes the BTM associated to the BSTM computed by $\P$
(\textit{cf}. Figure \ref{figABRMa}).

Let us denote by $S(t)$ the generating series of these trees counted 
by size. The generating 
serie statisfies the following functional equation (see
\href{http://oeis.org/A002212}{\texttt{A002212}} of \cite{oeis}):
\begin{eqnarray}
    S(t) &=& 1 + \frac{S(t)^2}{1-t} \text{ and thus } S(t) =
    \frac{1-t - \sqrt{5t^2 - 6t + 1}}{2t} \label{eqnFormClose}\\
    &=& 1+t+3t^2 +10t^3 +36t^4 +137t^5 +543t^6 +2219t^7 + \ldots\nonumber
\end{eqnarray}
This structure is in bijection with \emph{binary unary tree} structure.
Here is the list of trees of size 0,1,2 and 3:
\[\cdot, \input{arbres/ensArbr}\]

With Lemma \ref{lemmeAlgebraQuotient} and Theorem 
\ref{thmGoodHopfAlgebra} we know that the quotient of $\wqsym(\Ag)$ by the
ta\"iga relations has a natural basis indexed by $tass(\Ag^*)/_{\equiv_t}$
identified by BTM. 
We call \pbtm (planar binary tree with multiplicities) that 
quotient.
More precisely, we consider the basis $(\MM_u)_u$ of \wqsym obtained 
by the Hopf polynomial realization $r_{tass}$. We denote by $(\QQ^m_t)_t$ 
the canonical projection by the map $\pi$ of $(\MM_u)_u$ in \pbtm such 
that $\pi(\MM_u) := \QQ^m_t$ if $t = \B(u)$.
The product and coproduct are given by some explicit algorithms. 
For brevity, we only give here some examples.
 
      \input{arbres/quotient/prod}
      \input{arbres/quotient/coprod}
      \input{arbres/quotient/dual}

%% file: arbres/ensArbr.tex
\vcenter{\hbox{\scalebox{.6}{\begin{tikzpicture}[every node/.style={draw,circle}]
    \node (0) at (1/1.4,1) {$\scriptstyle 1$};
\end{tikzpicture}}}}, 
\vcenter{\hbox{\scalebox{.6}{\begin{tikzpicture}[every node/.style={draw,circle}]
    \node (0) at (1/1.4,2) {$\scriptstyle 1$};
    \node (1) at (2/1.4,1) {$\scriptstyle 1$};
\path[red,ultra thick]
    (0) edge (1)
;\end{tikzpicture}}}}, \vcenter{\hbox{\scalebox{.6}{\begin{tikzpicture}[every node/.style={draw,circle}]
    \node (1) at (1/1.4,1) {$\scriptstyle 1$};
    \node (0) at (2/1.4,2) {$\scriptstyle 1$};
\path[red,ultra thick]
    (0) edge (1)
;\end{tikzpicture}}}}, \vcenter{\hbox{\scalebox{.6}{\begin{tikzpicture}[every node/.style={draw,circle}]
    \node (0) at (1/1.4,1) {$\scriptstyle 2$};
\path[red,ultra thick]
;\end{tikzpicture}}}},
\vcenter{\hbox{\scalebox{.6}{\begin{tikzpicture}[every node/.style={draw,circle}]
    \node (0) at (1/1.4,3) {$\scriptstyle 1$};
    \node (1) at (2/1.4,2) {$\scriptstyle 1$};
    \node (2) at (3/1.4,1) {$\scriptstyle 1$};
\path[red,ultra thick]
    (1) edge (2)
    (0) edge (1)
;\end{tikzpicture}}}}, \vcenter{\hbox{\scalebox{.6}{\begin{tikzpicture}[every node/.style={draw,circle}]
    \node (0) at (1/1.4,3) {$\scriptstyle 1$};
    \node (2) at (2/1.4,1) {$\scriptstyle 1$};
    \node (1) at (3/1.4,2) {$\scriptstyle 1$};
\path[red,ultra thick]
    (1) edge (2)
    (0) edge (1)
;\end{tikzpicture}}}}, \vcenter{\hbox{\scalebox{.6}{\begin{tikzpicture}[every node/.style={draw,circle}]
    \node (0) at (1/1.4,2) {$\scriptstyle 1$};
    \node (1) at (2/1.4,1) {$\scriptstyle 2$};
\path[red,ultra thick]
    (0) edge (1)
;\end{tikzpicture}}}}, \vcenter{\hbox{\scalebox{.6}{\begin{tikzpicture}[every node/.style={draw,circle}]
    \node (1) at (1/1.4,1) {$\scriptstyle 1$};
    \node (0) at (2/1.4,2) {$\scriptstyle 1$};
    \node (2) at (3/1.4,1) {$\scriptstyle 1$};
\path[red,ultra thick]
    (0) edge (1)
    (0) edge (2)
;\end{tikzpicture}}}}, \vcenter{\hbox{\scalebox{.6}{\begin{tikzpicture}[every node/.style={draw,circle}]
    \node (1) at (1/1.4,2) {$\scriptstyle 1$};
    \node (2) at (2/1.4,1) {$\scriptstyle 1$};
    \node (0) at (3/1.4,3) {$\scriptstyle 1$};
\path[red,ultra thick]
    (1) edge (2)
    (0) edge (1)
;\end{tikzpicture}}}}, \vcenter{\hbox{\scalebox{.6}{\begin{tikzpicture}[every node/.style={draw,circle}]
    \node (2) at (1/1.4,1) {$\scriptstyle 1$};
    \node (1) at (2/1.4,2) {$\scriptstyle 1$};
    \node (0) at (3/1.4,3) {$\scriptstyle 1$};
\path[red,ultra thick]
    (1) edge (2)
    (0) edge (1)
;\end{tikzpicture}}}}, \vcenter{\hbox{\scalebox{.6}{\begin{tikzpicture}[every node/.style={draw,circle}]
    \node (1) at (1/1.4,1) {$\scriptstyle 2$};
    \node (0) at (2/1.4,2) {$\scriptstyle 1$};
\path[red,ultra thick]
    (0) edge (1)
;\end{tikzpicture}}}}, \vcenter{\hbox{\scalebox{.6}{\begin{tikzpicture}[every node/.style={draw,circle}]
    \node (0) at (1/1.4,2) {$\scriptstyle 2$};
    \node (1) at (2/1.4,1) {$\scriptstyle 1$};
\path[red,ultra thick]
    (0) edge (1)
;\end{tikzpicture}}}}, \vcenter{\hbox{\scalebox{.6}{\begin{tikzpicture}[every node/.style={draw,circle}]
    \node (1) at (1/1.4,1) {$\scriptstyle 1$};
    \node (0) at (2/1.4,2) {$\scriptstyle 2$};
\path[red,ultra thick]
    (0) edge (1)
;\end{tikzpicture}}}}, \vcenter{\hbox{\scalebox{.6}{\begin{tikzpicture}[every node/.style={draw,circle}]
    \node (0) at (1/1.4,1) {$\scriptstyle 3$};
\path[red,ultra thick]
;\end{tikzpicture}}}}

%% file: arbres/quotient/prod.tex
The product on the $(\QQ^m_t)_t$ basis is described thanks to the 
projection $\pi$. For example,
    \renewcommand{\QQm}[1]{\ensuremath{\QQ^m_{\vcenter{\hbox{\scalebox{.4}{\input{arbres/quotient/prod/#1}}}}}}}
    \definecolor{c1}{RGB}{0,169,157}
    \definecolor{c2}{RGB}{192,79,22}
    \definecolor{c3}{RGB}{148,0,212}
    \begin{eqnarray*}
        \pi(\mathbb{M}_{1312} \times \mathbb{M}_{1}) &=& \pi(
        \mathbb{M}_{13121} + \mathbb{M}_{13122} + \mathbb{M}_{13123} + 
        \mathbb{M}_{13124} + \mathbb{M}_{14123} + \mathbb{M}_{14132} + 
        \mathbb{M}_{24231})\\
        = \mathbb{Q}^m_{\vcenter{\hbox{\scalebox{.4}
{ { \newcommand{\nodea}{\node[draw,circle] (a) {$1$}
;}\newcommand{\nodeb}{\node[draw,circle] (b) {$2$}
;}\newcommand{\nodec}{\node[draw,circle] (c) {$1$}
;}\begin{tikzpicture}[auto]
\matrix[column sep=.1cm, row sep=.1cm,ampersand replacement=\&]{
         \& \nodea  \&         \\ 
 \nodeb  \&         \& \nodec  \\
};

\path[ultra thick, red] (a) edge (b) edge (c);
\end{tikzpicture}} }}}} \times \mathbb{Q}^m_{\vcenter{\hbox{\scalebox{.4}
{ { \newcommand{\nodea}{\node[draw,circle] (a) {$1$}
;}\begin{tikzpicture}[auto]
\matrix[column sep=.1cm, row sep=.1cm,ampersand replacement=\&]{
 \nodea  \\
};
\end{tikzpicture}} }}}} &=&
        \mathbb{Q}^m_{\vcenter{\hbox{\scalebox{.4}
{ { \newcommand{\nodea}{\node[draw,circle] (a) {$3$}
;}\newcommand{\nodeb}{\node[draw,circle] (b) {$1$}
;}\newcommand{\nodec}{\node[draw,circle] (c) {$1$}
;}\begin{tikzpicture}[auto]
\matrix[column sep=.1cm, row sep=.1cm,ampersand replacement=\&]{
 \nodea  \&         \&         \\ 
         \& \nodeb  \&         \\ 
         \&         \& \nodec  \\
};

\path[ultra thick, red] (b) edge (c)
    (a) edge (b);
\end{tikzpicture}} }}}} + \mathbb{Q}^m_{\vcenter{\hbox{\scalebox{.4}
{ { \newcommand{\nodea}{\node[draw,circle] (a) {$1$}
;}\newcommand{\nodeb}{\node[draw,circle] (b) {$1$}
;}\newcommand{\nodec}{\node[draw,circle] (c) {$2$}
;}\newcommand{\noded}{\node[draw,circle] (d) {$1$}
;}\begin{tikzpicture}[auto]
\matrix[column sep=.1cm, row sep=.1cm,ampersand replacement=\&]{
 \nodea  \&         \&         \\ 
         \& \nodeb  \&         \\ 
 \nodec  \&         \& \noded  \\
};

\path[ultra thick, red] (b) edge (c) edge (d)
    (a) edge (b);
\end{tikzpicture}} }}}} + \mathbb{Q}^m_{\vcenter{\hbox{\scalebox{.4}
{ { \newcommand{\nodea}{\node[draw,circle] (a) {$2$}
;}\newcommand{\nodeb}{\node[draw,circle] (b) {$2$}
;}\newcommand{\nodec}{\node[draw,circle] (c) {$1$}
;}\begin{tikzpicture}[auto]
\matrix[column sep=.1cm, row sep=.1cm,ampersand replacement=\&]{
         \& \nodea  \&         \\ 
 \nodeb  \&         \& \nodec  \\
};

\path[ultra thick, red] (a) edge (b) edge (c);
\end{tikzpicture}} }}}} + \mathbb{Q}^m_{\vcenter{\hbox{\scalebox{.4}
{ { \newcommand{\nodea}{\node[draw,circle] (a) {$1$}
;}\newcommand{\nodeb}{\node[draw,circle] (b) {$2$}
;}\newcommand{\nodec}{\node[draw,circle] (c) {$1$}
;}\newcommand{\noded}{\node[draw,circle] (d) {$1$}
;}\begin{tikzpicture}[auto]
\matrix[column sep=.1cm, row sep=.1cm,ampersand replacement=\&]{
         \& \nodea  \&         \&         \\ 
 \nodeb  \&         \& \nodec  \&         \\ 
         \&         \&         \& \noded  \\
};

\path[ultra thick, red] (c) edge (d)
    (a) edge (b) edge (c);
\end{tikzpicture}} }}}} + \mathbb{Q}^m_{\vcenter{\hbox{\scalebox{.4}
{ { \newcommand{\nodea}{\node[draw,circle] (a) {$2$}
;}\newcommand{\nodeb}{\node[draw,circle] (b) {$1$}
;}\newcommand{\nodec}{\node[draw,circle] (c) {$2$}
;}\begin{tikzpicture}[auto]
\matrix[column sep=.1cm, row sep=.1cm,ampersand replacement=\&]{
         \&         \& \nodea  \\ 
         \& \nodeb  \&         \\ 
 \nodec  \&         \&         \\
};

\path[ultra thick, red] (b) edge (c)
    (a) edge (b);
\end{tikzpicture}} }}}} + \mathbb{Q}^m_{\vcenter{\hbox{\scalebox{.4}
{ { \newcommand{\nodea}{\node[draw,circle] (a) {$1$}
;}\newcommand{\nodeb}{\node[draw,circle] (b) {$1$}
;}\newcommand{\nodec}{\node[draw,circle] (c) {$2$}
;}\newcommand{\noded}{\node[draw,circle] (d) {$1$}
;}\begin{tikzpicture}[auto]
\matrix[column sep=.1cm, row sep=.1cm,ampersand replacement=\&]{
         \&         \& \nodea  \&         \\ 
         \& \nodeb  \&         \& \noded  \\ 
 \nodec  \&         \&         \&         \\
};

\path[ultra thick, red] (b) edge (c)
    (a) edge (b) edge (d);
\end{tikzpicture}} }}}} + \mathbb{Q}^m_{\vcenter{\hbox{\scalebox{.4}
{ { \newcommand{\nodea}{\node[draw,circle] (a) {$1$}
;}\newcommand{\nodeb}{\node[draw,circle] (b) {$1$}
;}\newcommand{\nodec}{\node[draw,circle] (c) {$2$}
;}\newcommand{\noded}{\node[draw,circle] (d) {$1$}
;}\begin{tikzpicture}[auto]
\matrix[column sep=.1cm, row sep=.1cm,ampersand replacement=\&]{
         \&         \& \nodea  \\ 
         \& \nodeb  \&         \\ 
 \nodec  \&         \& \noded  \\
};

\path[ultra thick, red] (b) edge (c) edge (d)
    (a) edge (b);
\end{tikzpicture}} }}}}\;.
    \end{eqnarray*}
%

%% file: arbres/quotient/coprod.tex
Similary, the coproduct is described thanks to the projection $\pi$. We could 
remark that it is the same as the coproduct on $(\QQ_t)_t$ of \pbt keeping the 
multiplicities with nodes. For example,  
    \begin{eqnarray*}
        \pi(\Delta(\mathbb{M}_{3112})) &=& \pi(
            1 \otimes \mathbb{M}_{3112} + 
            \mathbb{M}_{11} \otimes \mathbb{M}_{21} + 
            \mathbb{M}_{112} \otimes \mathbb{M}_{1} + 
            \mathbb{M}_{3112} \otimes 1
           )\\ 
        = \Delta(\mathbb{Q}^m_{\vcenter{\hbox{\scalebox{.4}
{ { \newcommand{\nodea}{\node[draw,circle] (a) {$1$}
;}\newcommand{\nodeb}{\node[draw,circle] (b) {$2$}
;}\newcommand{\nodec}{\node[draw,circle] (c) {$1$}
;}\begin{tikzpicture}[auto]
\matrix[column sep=.3cm, row sep=.3cm,ampersand replacement=\&]{
         \& \nodea  \&         \\ 
 \nodeb  \&         \& \nodec  \\
};

\path[ultra thick, red] (a) edge (b) edge (c);
\end{tikzpicture}} }}}}) &=& 
 1 \otimes \mathbb{Q}^m_{\vcenter{\hbox{\scalebox{.4}
{ { \newcommand{\nodea}{\node[draw,circle] (a) {$1$}
;}\newcommand{\nodeb}{\node[draw,circle] (b) {$2$}
;}\newcommand{\nodec}{\node[draw,circle] (c) {$1$}
;}\begin{tikzpicture}[auto]
\matrix[column sep=.3cm, row sep=.3cm,ampersand replacement=\&]{
         \& \nodea  \&         \\ 
 \nodeb  \&         \& \nodec  \\
};

\path[ultra thick, red] (a) edge (b) edge (c);
\end{tikzpicture}} }}}} + \mathbb{Q}^m_{\vcenter{\hbox{\scalebox{.4}
{ { \newcommand{\nodea}{\node[draw,circle] (a) {$2$}
;}\begin{tikzpicture}[auto]
\matrix[column sep=.3cm, row sep=.3cm,ampersand replacement=\&]{
 \nodea  \\
};
\end{tikzpicture}} }}}} \otimes \mathbb{Q}^m_{\vcenter{\hbox{\scalebox{.4}
{ { \newcommand{\nodea}{\node[draw,circle] (a) {$1$}
;}\newcommand{\nodeb}{\node[draw,circle] (b) {$1$}
;}\begin{tikzpicture}[auto]
\matrix[column sep=.3cm, row sep=.3cm,ampersand replacement=\&]{
 \nodea  \&         \\ 
         \& \nodeb  \\
};

\path[ultra thick, red] (a) edge (b);
\end{tikzpicture}} }}}} + \mathbb{Q}^m_{\vcenter{\hbox{\scalebox{.4}
{ { \newcommand{\nodea}{\node[draw,circle] (a) {$1$}
;}\newcommand{\nodeb}{\node[draw,circle] (b) {$2$}
;}\begin{tikzpicture}[auto]
\matrix[column sep=.3cm, row sep=.3cm,ampersand replacement=\&]{
         \& \nodea  \\ 
 \nodeb  \&         \\
};

\path[ultra thick, red] (a) edge (b);
\end{tikzpicture}} }}}} \otimes \mathbb{Q}^m_{\vcenter{\hbox{\scalebox{.4}
{ { \newcommand{\nodea}{\node[draw,circle] (a) {$1$}
;}\begin{tikzpicture}[auto]
\matrix[column sep=.3cm, row sep=.3cm,ampersand replacement=\&]{
 \nodea  \\
};
\end{tikzpicture}} }}}} + \mathbb{Q}^m_{\vcenter{\hbox{\scalebox{.4}
{ { \newcommand{\nodea}{\node[draw,circle] (a) {$1$}
;}\newcommand{\nodeb}{\node[draw,circle] (b) {$2$}
;}\newcommand{\nodec}{\node[draw,circle] (c) {$1$}
;}\begin{tikzpicture}[auto]
\matrix[column sep=.3cm, row sep=.3cm,ampersand replacement=\&]{
         \& \nodea  \&         \\ 
 \nodeb  \&         \& \nodec  \\
};

\path[ultra thick, red] (a) edge (b) edge (c);
\end{tikzpicture}} }}}} \otimes 1\;.
    \end{eqnarray*}

%% file: arbres/quotient/dual.tex
We consider $\pbtm^\#$ the dual of \pbtm. This is a sub-algebra
of $\wqsym^\#$. We denote by $(\PP^m_t)_t := (\QQ^m_t)^\#$ its dual basis:
$\langle \QQ^m_t, \PP^m_{t'}\rangle = \delta_{t,t'}$.
The product is given by: 
\begin{eqnarray}
    \PP_{t'}^m \times \PP_{t''}^m &=& \sum_t \langle \Delta(\QQ^m_t),
    \PP_{t'}^m \times \PP_{t''}^m\rangle \PP_t^m\;.
\end{eqnarray} 
Here is an example,
    \begin{eqnarray*}
        \mathbb{P}^m_{\vcenter{\hbox{\scalebox{.4}
{ { \newcommand{\nodea}{\node[draw,circle] (a) {$3$}
;}\newcommand{\nodeb}{\node[draw,circle] (b) {$1$}
;}\begin{tikzpicture}[auto]
\matrix[column sep=.1cm, row sep=.1cm,ampersand replacement=\&]{
 \nodea  \&         \\ 
         \& \nodeb  \\
};

\path[ultra thick, red] (a) edge (b);
\end{tikzpicture}} }}}} \times \mathbb{P}^m_{\vcenter{\hbox{\scalebox{.4}
{ { \newcommand{\nodea}{\node[draw,circle] (a) {$4$}
;}\newcommand{\nodeb}{\node[draw,circle] (b) {$2$}
;}\begin{tikzpicture}[auto]
\matrix[column sep=.1cm, row sep=.1cm,ampersand replacement=\&]{
         \& \nodea \\ 
 \nodeb  \\
};

\path[ultra thick, red] (a) edge (b);
\end{tikzpicture}} }}}} &=& \mathbb{P}^m_{\vcenter{\hbox{\scalebox{.4}
{ { \newcommand{\nodea}{\node[draw,circle] (a) {$3$}
;}\newcommand{\nodeb}{\node[draw,circle] (b) {$1$}
;}\newcommand{\nodec}{\node[draw,circle] (c) {$4$}
;}\newcommand{\noded}{\node[draw,circle] (d) {$2$}
;}\begin{tikzpicture}[auto]
\matrix[column sep=.1cm, row sep=.1cm,ampersand replacement=\&]{
 \nodea  \\ 
         \& \nodeb \\ 
         \&         \& \nodec  \\ 
         \& \noded  \\
};

\path[ultra thick, red] (c) edge (d)
    (b) edge (c)
    (a) edge (b);
\end{tikzpicture}} }}}} + \mathbb{P}^m_{\vcenter{\hbox{\scalebox{.4}
{ { \newcommand{\nodea}{\node[draw,circle] (a) {$3$}
;}\newcommand{\nodeb}{\node[draw,circle] (b) {$4$}
;}\newcommand{\nodec}{\node[draw,circle] (c) {$1$}
;}\newcommand{\noded}{\node[draw,circle] (d) {$2$}
;}\begin{tikzpicture}[auto]
\matrix[column sep=.1cm, row sep=.1cm,ampersand replacement=\&]{
 \nodea  \&         \\ 
         \& \nodeb  \\ 
 \nodec  \&         \\ 
         \& \noded  \\
};

\path[ultra thick, red] (c) edge (d)
    (b) edge (c)
    (a) edge (b);
\end{tikzpicture}} }}}} + \mathbb{P}^m_{\vcenter{\hbox{\scalebox{.4}
{ { \newcommand{\nodea}{\node[draw,circle] (a) {$3$}
;}\newcommand{\nodeb}{\node[draw,circle] (b) {$4$}
;}\newcommand{\nodec}{\node[draw,circle] (c) {$2$}
;}\newcommand{\noded}{\node[draw,circle] (d) {$1$}
;}\begin{tikzpicture}[auto]
\matrix[column sep=.1cm, row sep=.1cm,ampersand replacement=\&]{
         \& \nodea  \&         \\ 
         \&         \& \nodeb  \\ 
         \& \nodec  \&         \\ 
 \noded  \&         \&         \\
};

\path[ultra thick, red] (c) edge (d)
    (b) edge (c)
    (a) edge (b);
\end{tikzpicture}} }}}} + \mathbb{P}^m_{\vcenter{\hbox{\scalebox{.4}
{ { \newcommand{\nodea}{\node[draw,circle] (a) {$4$}
;}\newcommand{\nodeb}{\node[draw,circle] (b) {$3$}
;}\newcommand{\nodec}{\node[draw,circle] (c) {$1$}
;}\newcommand{\noded}{\node[draw,circle] (d) {$2$}
;}\begin{tikzpicture}[auto]
\matrix[column sep=.1cm, row sep=.1cm,ampersand replacement=\&]{
         \& \nodea  \\ 
 \nodeb  \\ 
         \& \nodec  \\ 
         \&         \& \noded \\
};

\path[ultra thick, red] (c) edge (d)
    (b) edge (c)
    (a) edge (b);
\end{tikzpicture}} }}}} + \mathbb{P}^m_{\vcenter{\hbox{\scalebox{.4}
{ { \newcommand{\nodea}{\node[draw,circle] (a) {$4$}
;}\newcommand{\nodeb}{\node[draw,circle] (b) {$3$}
;}\newcommand{\nodec}{\node[draw,circle] (c) {$2$}
;}\newcommand{\noded}{\node[draw,circle] (d) {$1$}
;}\begin{tikzpicture}[auto]
\matrix[column sep=.1cm, row sep=.1cm,ampersand replacement=\&]{
         \& \nodea \\ 
 \nodeb  \&        \\ 
         \& \nodec \\ 
 \noded  \&        \\
};

\path[ultra thick, red] (c) edge (d)
    (b) edge (c)
    (a) edge (b);
\end{tikzpicture}} }}}} + \mathbb{P}^m_{\vcenter{\hbox{\scalebox{.4}
{ { \newcommand{\nodea}{\node[draw,circle] (a) {$4$}
;}\newcommand{\nodeb}{\node[draw,circle] (b) {$2$}
;}\newcommand{\nodec}{\node[draw,circle] (c) {$3$}
;}\newcommand{\noded}{\node[draw,circle] (d) {$1$}
;}\begin{tikzpicture}[auto]
\matrix[column sep=.1cm, row sep=.1cm,ampersand replacement=\&]{
         \&         \& \nodea  \\ 
         \& \nodeb  \&         \\ 
 \nodec  \&         \&         \\ 
         \& \noded  \&         \\
};

\path[ultra thick, red] (c) edge (d)
    (b) edge (c)
    (a) edge (b);
\end{tikzpicture}} }}}}\;.
    \end{eqnarray*}
If we consider only shape tree, the product is exactly 
the product of $(\PP_t)_t$ basis in \pbt \cite{hivert2002analogue}. Hence
this product is a shifted shuffle on trees.\\

The coproduct is given by:
    $\Delta^\#(\PP_t^m) = \sum_{t',t''}\ \langle \QQ^m_{t'} \times
    \QQ^m_{t''}, \PP_t^m \rangle\ \PP^m_{t'} \otimes \PP^m_{t''}$.
Here is an example:
    
    \begin{eqnarray*}
        \Delta^\#(\mathbb{P}^m_{\vcenter{\hbox{\scalebox{.4}
{ { \newcommand{\nodea}{\node[draw,circle] (a) {$2$}
;}\newcommand{\nodeb}{\node[draw,circle] (b) {$1$}
;}\begin{tikzpicture}[auto]
\matrix[column sep=.1cm, row sep=.1cm,ampersand replacement=\&]{
         \& \nodea  \\ 
 \nodeb  \&         \\
};

\path[ultra thick, red] (a) edge (b);
\end{tikzpicture}} }}}}) &=& 1 \otimes \mathbb{P}^m_{\vcenter{\hbox{\scalebox{.4}
{ { \newcommand{\nodea}{\node[draw,circle] (a) {$2$}
;}\newcommand{\nodeb}{\node[draw,circle] (b) {$1$}
;}\begin{tikzpicture}[auto]
\matrix[column sep=.1cm, row sep=.1cm,ampersand replacement=\&]{
         \& \nodea  \\ 
 \nodeb  \&         \\
};

\path[ultra thick, red] (a) edge (b);
\end{tikzpicture}} }}}} + \mathbb{P}^m_{\vcenter{\hbox{\scalebox{.4}
{ { \newcommand{\nodea}{\node[draw,circle] (a) {$1$}
;}\begin{tikzpicture}[auto]
\matrix[column sep=.1cm, row sep=.1cm,ampersand replacement=\&]{
 \nodea  \\
};
\end{tikzpicture}} }}}}
        \otimes \bigg(\mathbb{P}^m_{\vcenter{\hbox{\scalebox{.4}
{ { \newcommand{\nodea}{\node[draw,circle] (a) {$2$}
;}\begin{tikzpicture}[auto]
\matrix[column sep=.1cm, row sep=.1cm,ampersand replacement=\&]{
 \nodea  \\
};
\end{tikzpicture}} }}}} + \mathbb{P}^m_{\vcenter{\hbox{\scalebox{.4}
{ { \newcommand{\nodea}{\node[draw,circle] (a) {$1$}
;}\newcommand{\nodeb}{\node[draw,circle] (b) {$1$}
;}\begin{tikzpicture}[auto]
\matrix[column sep=.1cm, row sep=.1cm,ampersand replacement=\&]{
         \& \nodea  \\ 
 \nodeb  \&         \\
};

\path[ultra thick, red] (a) edge (b);
\end{tikzpicture}} }}}}\bigg) + 
        \bigg(\mathbb{P}^m_{\vcenter{\hbox{\scalebox{.4}
{ { \newcommand{\nodea}{\node[draw,circle] (a) {$1$}
;}\newcommand{\nodeb}{\node[draw,circle] (b) {$1$}
;}\begin{tikzpicture}[auto]
\matrix[column sep=.1cm, row sep=.1cm,ampersand replacement=\&]{
 \nodea  \&         \\ 
         \& \nodeb  \\
};

\path[ultra thick, red] (a) edge (b);
\end{tikzpicture}} }}}} + \mathbb{P}^m_{\vcenter{\hbox{\scalebox{.4}
{ { \newcommand{\nodea}{\node[draw,circle] (a) {$1$}
;}\newcommand{\nodeb}{\node[draw,circle] (b) {$1$}
;}\begin{tikzpicture}[auto]
\matrix[column sep=.1cm, row sep=.1cm,ampersand replacement=\&]{
         \& \nodea  \\ 
 \nodeb  \&         \\
};

\path[ultra thick, red] (a) edge (b);
\end{tikzpicture}} }}}}\bigg) \otimes \mathbb{P}^m_{\vcenter{\hbox{\scalebox{.4}
{ { \newcommand{\nodea}{\node[draw,circle] (a) {$1$}
;}\begin{tikzpicture}[auto]
\matrix[column sep=.1cm, row sep=.1cm,ampersand replacement=\&]{
 \nodea  \\
};
\end{tikzpicture}} }}}}
        + \mathbb{P}^m_{\vcenter{\hbox{\scalebox{.4}
{ { \newcommand{\nodea}{\node[draw,circle] (a) {$2$}
;}\newcommand{\nodeb}{\node[draw,circle] (b) {$1$}
;}\begin{tikzpicture}[auto]
\matrix[column sep=.1cm, row sep=.1cm,ampersand replacement=\&]{
         \& \nodea  \\ 
 \nodeb  \&         \\
};

\path[ultra thick, red] (a) edge (b);
\end{tikzpicture}} }}}} \otimes 1\;.
    \end{eqnarray*}
   

%% file: wip/hook.tex
Its well known from \cite{knuth1970permutations} (\textsection 
5.14 ex. 20) that the number of decreasing labelling of a binary tree 
is given by a simple product formula. \cite{hivert2002analogue} remarks 
that this is also the number of permutations given upon a tree by the binary 
search tree insertion. In this section we generalize this formula for trees
with multiplicities.

\begin{prop}
    The cardinal $f(T)$ of the ta\"iga class associated to $T$
    (\textit{i.e.} the set of packed words giving the tree $T$ by the 
    insertion algorithm $\B$) is 
    given by 
    \begin{eqnarray}
        f(T) &=& \vabs T!\left(\prod_{t\in T} \vabs t (m(t)-1)!\right)^{-1}\;.
        \label{eqnFormEq}
    \end{eqnarray}
    where $t$ ranges throwgh all the subtrees of $T$ and 
    $\vabs T$ denotes the size of $T$ (the sum of the multiplicities).
\end{prop}

\begin{exemple}
    The ta\"iga class of $T:=\vcenter{\hbox{\scalebox{.6}{{ \newcommand{\nodea}{\node[draw,circle] (a) {$2$}
;}\newcommand{\nodeb}{\node[draw,circle] (b) {$1$}
;}\newcommand{\nodec}{\node[draw,circle] (c) {$2$}
;}\begin{tikzpicture}[auto]
\matrix[column sep=.3cm, row sep=.3cm,ampersand replacement=\&]{
         \& \nodea  \&         \\ 
 \nodeb  \&         \& \nodec  \\
};

\path[ultra thick, red] (a) edge (b) edge (c);
\end{tikzpicture}}}}}$ contains 
    $12$ packed words $w$:
    \begin{eqnarray*}
        {23132, 33122, 31232, 32312, 13232, 33212, 23312, 32132, 21332, 31322, 
        12332, 13322}\;.
    \end{eqnarray*}
    The class of $\vcenter{\hbox{\scalebox{.6}{{ \newcommand{\nodea}{\node[draw,circle] (a) {$2$}
;}\newcommand{\nodeb}{\node[draw,circle] (b) {$2$}
;}\newcommand{\nodec}{\node[draw,circle] (c) {$7$}
;}\newcommand{\noded}{\node[draw,circle] (d) {$1$}
;}\newcommand{\nodee}{\node[draw,circle] (e) {$4$}
;}\newcommand{\nodef}{\node[draw,circle] (f) {$2$}
;}\begin{tikzpicture}[auto]
\matrix[column sep=.3cm, row sep=.3cm,ampersand replacement=\&]{
         \&         \&         \& \nodea  \&         \&         \&         \\ 
         \& \nodeb  \&         \&         \&         \& \noded  \&         \\ 
 \nodec  \&         \&         \&         \& \nodee  \&         \& \nodef  \\
};

\path[ultra thick, red] (b) edge (c)
    (d) edge (e) edge (f)
    (a) edge (b) edge (d);
\end{tikzpicture}}}}}$ 
    contains $23,337,600 = \frac{18!}{\left(18 \cdot
        9 \cdot 7 \cdot 7 \cdot 4 \cdot 2 \right) \left(
        1! 1! 6! 0! 3! 1!\right)}$ packed words.
\end{exemple}

This formula is easily proven by induction. However, we prefer to give
a generating series proof as in \cite{hivert2008trees}.
Let $\A$ be an associative algebra, and consider the functional equation
for power series $x\in \A[[z]]$:
\begin{eqnarray}
    x &=& a + \sum_{k\geqslant 1} B_k(x,x)\;, \label{eqnFuncFeq}
\end{eqnarray}
where $a \in \A$ and for any $k>0$, $B_k(x,y)$ is a bilinear map with
values in $\A[[z]]$. We suppose such that the valuation of $B_k(x,y)$ is strictly
greater than the sum of the valuations of $x$ and $y$ (plus $k$). 
Then Equation \ref{eqnFuncFeq} has a unique solution:
\begin{eqnarray}
    x &=& a + \sum_{k\geqslant 1} \bigg(B_k(a,a) + B_k(a, \sum_{k'\geqslant 1} B_{k'}(a,a)) +
    B_k(\sum_{k'\geqslant 1} B_{k'}(a,a)), a) + \ldots \bigg) \label{eqnFuncSolu}\\
    &=& \sum_{T\in \textbf{BTM}} B_T(a)\;, \nonumber 
\end{eqnarray}
where for a tree $T$, $B_T(a)$ is the result of evaluating the expression
formed by labelling by $a$ the leaves of the complete tree associated to
$T$ and by $B_k$ its internal node labelled by $k$.

\noindent
For example:
    $B_{\vcenter{\hbox{\scalebox{.6}{{ \newcommand{\nodea}{\node[draw,circle] (a) {$3$}
;}\newcommand{\nodeb}{\node[draw,circle] (b) {$6$}
;}\newcommand{\nodec}{\node[draw,circle] (c) {$2$}
;}\newcommand{\noded}{\node[draw,circle] (d) {$2$}
;}\begin{tikzpicture}[auto]
\matrix[column sep=.3cm, row sep=.3cm,ampersand replacement=\&]{
         \&         \&         \& \nodea  \&         \\ 
         \& \nodeb  \&         \&         \& \noded  \\ 
         \&         \& \nodec  \&         \&         \\
};

\path[ultra thick, red] (b) edge (c)
    (a) edge (b) edge (d);
\end{tikzpicture}}}}}}(a) 
    = B_3(B_6(a, B_2(a,a)), B_2(a,a))$.

So if we try to solve the fixed point problem:
\begin{eqnarray}
    x = 1 + \int_{0}^{z} e^s x(s)^2 ds =  
    1 + \sum_{k\geqslant 1} \int_0^z \frac{s^{k-1}}{(k-1)!}x(s)^2ds
    &=& 1 + \sum_{k \geqslant 1} B_k(x,x)\;, \label{eqnarrayFixedPointProblem} 
\end{eqnarray}
where $B_k(x,y) = \int_0^z \frac{s^{k-1}}{(k-1)!}x(s)y(s)ds$.
Then for a binary tree of non-negative integer $T$, $B_T(1)$ is the 
monomial obtained by putting $1$ on each leaf and integrating at each 
node $n$ the product of the evaluations of its subtrees and $s^k/k!$ 
with $m(n) = k+1$.

For example: $\vcenter{\hbox{\scalebox{.7}{\input{wip/ex}}}}$.
    
One can observe that $B_T(1) = f(T) \frac{z^n}{n!}$, where $n = \vabs T$.

To prove the \emph{hook length formula}, following the same technique
as in \cite{hivert2008trees}, we want to lift in $\wqsym^\#$ the 
fixed point computation of Equation \ref{eqnarrayFixedPointProblem}.
Recall from \cite{hivert1999combinatoire, novelli2005construction} that the
multiplication rule of the dual basis $\SS_u$ ($\MM_u^\# := \SS_u$) is,
for $u, v$ two packed words of respectively size $k$ and $l$,
$\SS_u \cdot \SS_v = \sumOld_{w \in u \cshuffle v} \SS_w$
with $u \cshuffle v$ the set of words apparing in the shifted shuffle 
of $u$ and $v$. 
The number of terms is a binomial: $\binom{k+l}{l}$ for $k$ and 
$l$ the length of $u$ and $v$.
Hence, the linear map $\varphi : \SS_u \mapsto \frac{z^n}{n!}$
with $n$ the length of $u$ is a morphism of algebras from $\wqsym^\#$ to 
$\KK[[z]]$. 
%
For $u,v$ two packed words of respective size $n-1$ and $m$, set
\begin{eqnarray}
   B_k(\SS_u, \SS_v) &:=& \sum_{w \in (u \cshuffle 1^{k-1}\cshuffle 
   v)\cdot n} \SS_w\;.
\end{eqnarray} 
The crucial observation which allows to express the \emph{hook length
formula} in a generating series way is the following theorem:
\begin{thm}
    \label{thmHookForm}
    In the binary tree (with multiplicity) solution (Equation 
    \ref{eqnFuncSolu})
    of Equation \ref{eqnarrayFixedPointProblem},
    \begin{eqnarray}
        B_T(1) &=& \sum_{\B(u) = T} \SS_u,
    \end{eqnarray}
    In particular, $B_T(1)$ coincide with $\PP^m_T$, the natural basis 
    of $\pbtm^\#$.
\end{thm}

\begin{coro}
    The number of packed words $u$ such that $\B(u) = T$ is computed by 
    $f(T)$.  
\end{coro}

%% file: wip/ex.tex
{\newcommand{\AOnooo}{\node (0) [circle, fill=white, text=black, draw] 
    {$3$};}\newcommand{\AOnooa}{\node (1) [circle, fill=white, text=black, draw] 
    {$2$};}\newcommand{\AOnoob}{\node (2) [circle, fill=white, text=black, draw] 
    {$2$};}\newcommand{\AOnooc}{\node (3) [circle, fill=white, text=black, draw] 
    {$1$};}\newcommand{\AOnood}{\node (4) [circle, fill=white, text=black, draw] 
    {$1$};}\newcommand{\AOnooe}{\node (5) [circle, fill=white, text=black, draw] 
    {$1$};}\newcommand{\AOnoof}{\node (6) 
    {$\frac{z^7}{16\times7}$};}\newcommand{\AOnoog}{\node (7)  
    {$\frac{z^4}{8}$};}\newcommand{\AOnooh}{\node (8)  
    {$\frac{z^2}{2}$};}\newcommand{\AOnooi}{\node (9)  
    {$z$};}\newcommand{\AOnooj}{\node (10) 
    {$\frac{z^2}{2}$};}\newcommand{\AOnook}{\node (11) 
    {$z$};}\newcommand{\AOno}{\node (12) 
    {$1$};}\newcommand{\AOna}{\node (13) {$1$};}\newcommand{\AOnb}{\node (14) {$1$};}\newcommand{\AOnc}{\node (15) {$1$};}\newcommand{\AOnd}{\node (16) {$1$};}\newcommand{\AOne}{\node (17) {$1$};}\newcommand{\AOnf}{\node (18) {$1$};}\begin{tikzpicture}[auto,scale=.3]
\matrix[column sep=.1cm, row sep=.1cm]{
         &         & \AOnooo &         &         &&                  &&         &       &         & \AOnoof \\ 
         & \AOnooa &         & \AOnood &&&\node{$\xrightarrow{B_T}$};&&         &       & \AOnoog &         &       & \AOnooh \\ 
 \AOnooc &         & \AOnoob &         & \AOnooe &&                  && \AOnooi &       &         & \AOnooj & \AOnf &         & \AOnook\\
         &         &         &         &         &&             &\AOno&         & \AOna &  \AOnb  &         & \AOnc &  \AOnd  &         & \AOne\\         
};
\path[red, ultra thick] (1) edge (3) edge (2)
                        (0) edge (4) edge (1)
                        (4) edge (5);
\path[blue, thick] (6) edge (7) edge (8)
                        (7) edge (9) edge (10)
                        (8) edge (11) edge (18)
                        (9) edge (12) edge (13)
                        (10)edge (14) edge (15)
                        (11)edge (16) edge (17);
\end{tikzpicture}}

%% file: perspectives.tex
In this paper, we unraveled some new combinatorics on binary trees with
multiplicities from the union of the \emph{sylvester} and \emph{stalactic} monoids. Using
the machinery of \emph{realizations}, we built a Hopf algebra on those trees,
allowing us to give a generating series proof of a new \emph{hook length
formula}. Following \cite{hivert2008trees}, it is very likely that we 
will also be able to prove a \emph{$q$-hook length formula}. On the other hand, 
the usual case of the
\nom{Loday-Ronco} algebra has a lot of nice properties. For example, the product and
coproduct can be expressed by the means of an order on the trees called the
\emph{Tamari Lattice} \cite{loday1998hopf}. It would be good to know if such 
a lattice exists for trees with multiplicities. This should also relate 
to N. \nom{Reading} work on
lattice congruences \cite{reading2005lattice}. Also it could be 
interesting to study some other
combinations in the lattice of \emph{good monoids}. For example, the union of the
\emph{plactic monoid} and the \emph{stalactic monoid} should give a 
Hopf algebra of
\emph{tableaux} with multiplicties. Finally, in our construction, it seems that $std$,
$tass$ and $park$ play some canonical role from which everything else is
built. Are there some more examples? Is there a definition for such a 
$\phi$-map? Could we except to always have a \emph{hook formula} as soon as
we have a \emph{good monoid}?